\documentclass{article}
\usepackage[utf8]{inputenc}
\usepackage[english]{babel}
\usepackage{tikz}
\usepackage{amsmath}
\usepackage{amssymb}
\usepackage{latexsym}
\usepackage{MnSymbol}
\usepackage[strict]{changepage}
\usepackage{tkz-euclide}
\usepackage{comment}
\usepackage{amsthm}
\newtheorem*{remark}{Remark}
\makeatletter
\renewcommand\@biblabel[1]{#1.}
\author{BENJAMIN EDUN}
\title{Finite and Infinite Nested Square Roots Convergent to Unity}
\date{}
\begin{document}
\maketitle
\begin{adjustwidth}{4em}{4em}
\textbf{\normalsize{Abstract}}.  We investigate finite and infinite nested square root formulas convergent to unity.
\end{adjustwidth}\vspace{3mm}
\textbf{\large{Keywords:}}Nested square roots, Half angle formula,  Recursive Formula, Binomial Expansion.\vspace{3mm}\newline
{\bf 2010 Mathematics Subject Classification}:11Y60 
\section{Introduction}
Most of the convergent nested square root formulas discussed in mathematics literature are in a form in which the real numbers within the nested square roots have integer exponents.  
As compiled in  {\cite{Servi}}, below are examples
\newline
\begin{align}
\pi &=\lim_{n \to \infty}\left[   2^{n} \underbrace{\sqrt{ 2-\sqrt{ 2+\cdots+\sqrt{2+\sqrt{2}}}}}_{n \hspace{1mm}square \hspace{1mm} roots} \right],\label{eq:1}\\
\pi &=\lim_{n \to \infty} \left[ 3 \cdot 2^{n-1}\underbrace{\sqrt{ 2-\sqrt{2+\cdots +\sqrt{ 2+\sqrt{3}}}}}_{n \hspace{1mm}square \hspace{1mm} roots} \right],\label{eq:2}\\
\pi &=\lim_{n \to \infty}  \left[\frac{4}{3}\cdot 2^{n-2} \underbrace{\sqrt{ 2-\sqrt{2+\cdots+\sqrt{ 2-\sqrt{2}}}}}_{n \hspace{1mm}square \hspace{1mm} roots}\right],\label{eq:3}\\
\pi &=\lim_{n \to \infty} \left[\frac{3}{5}\cdot  2^{n-1} \underbrace{\sqrt{ 2-\sqrt{2+\cdots+\sqrt{ 2-\sqrt{3}}}}}_{n \hspace{1mm}square \hspace{1mm} roots}\right]\label{eq:4},
\end{align}
It is quite clear that all numbers within the nested square roots in the expressions (\ref{eq:1}), (\ref{eq:2}), (\ref{eq:3}) and (\ref{eq:4}) have an exponent of unity. These expressions can be derived from the half angle formula,
\begin{equation*}
    sin(x)=\sqrt{\frac{1-cos(2x)}{2}}
\end{equation*}
Using the half-angle formula, 
\begin{equation}\label{eq:5}  
sin\left(\frac{\pi}{2^{n+1}}\right)= \left[\frac{1}{2} \underbrace{\sqrt{ 2-\sqrt{ 2+\cdots+\sqrt{2+\sqrt{2}}}}}_{n \hspace{1mm}square \hspace{1mm} roots} \right],
\end{equation} 
\begin{equation}\label{eq:6}
    \lim_{n\to \infty}sin\left(\frac{\pi}{2^{n+1}}\right)=\lim_{n \to \infty} \frac{\pi}{2^{n+1}}= \left[\frac{1}{2} \underbrace{\sqrt{ 2-\sqrt{ 2+\cdots+\sqrt{2+\sqrt{2}}}}}_{n \hspace{1mm}square \hspace{1mm} roots} \right],
\end{equation}
Thus,
\begin{equation}\label{eq:7}
 \pi=\lim_{n \to \infty}\left[   2^{n} \underbrace{\sqrt{ 2-\sqrt{ 2+\cdots+\sqrt{2+\sqrt{2}}}}}_{n \hspace{1mm}square \hspace{1mm} roots} \right]  
\end{equation}
Also from (\ref{eq:5}),
\begin{equation}\label{eq:8} 
sin\left(\frac{\arccos\frac{\sqrt{2}}{2}}{2^{n-1}}\right)= \left[\frac{1}{2} \underbrace{\sqrt{ 2-\sqrt{ 2+\cdots+\sqrt{2+\sqrt{2}}}}}_{n \hspace{1mm}square \hspace{1mm} roots} \right],
\end{equation}
\begin{equation}\label{eq:9} 
\lim_{n\to \infty} \frac{sin\left(\frac{\arccos\frac{\sqrt{2}}{2}}{2^{n-1}}\right)}{\left(\frac{\arccos\frac{\sqrt{2}}{2}}{2^{n-1}}\right)}=\lim_{n \to \infty} \left[\left(\frac{2^{n-2}}{\arccos\frac{\sqrt{2}}{2}}\right) \underbrace{\sqrt{ 2-\sqrt{ 2+\cdots+\sqrt{2+\sqrt{2}}}}}_{n \hspace{1mm}square \hspace{1mm} roots} \right], 
\end{equation}
Then, 
\begin{equation}\label{eq:10} 
\lim_{n \to \infty} \left[(\frac{2^{n-2}}{arcos\frac{\sqrt{2}}{2}}\Bigg) \underbrace{\sqrt{ 2-\sqrt{ 2+\cdots+\sqrt{2+\sqrt{2}}}}}_{n \hspace{1mm}square \hspace{1mm} roots} \right]=1,
\end{equation}
Similarly we could derive (\ref{eq:2}) , (\ref{eq:3}) and (\ref{eq:4}), and their corresponding forms equivalent to unity, using the half-angle formula. However, our focus is to derive infinitely nested square root formulas in the form
\begin{equation}\label{eq:11}
\lim_{n \to \infty}\left[g\bigg(\frac{s}{m}, n\bigg)
\underbrace{
\sqrt{ (2^{k_n}m)^{q_n}-\sqrt{(2^{k_{n-1}})^{q_{n-1}}+\cdots+\sqrt{m+s}}}}_{n \hspace{1mm}square \hspace{1mm} roots}
\right]=1
\end{equation}
and finite nested square root formulas in the form,
\begin{equation}\label{eq:12}
\lim_{m \to \infty}\left[g\left(\frac{s}{m}, n\right)
\underbrace{
\sqrt{ (2^{k_n}m)^{q_n}-\sqrt{(2^{k_{n-1}})^{q_{n-1}}+\cdots+\sqrt{m+s}}}}_{n \hspace{1mm}square \hspace{1mm} roots}
\right]=1,
\end{equation}
in which the numbers within the nested square roots can have fractional exponents.
\section{Infinitely Nested Square Roots Convergent to Unity}
Let the half angle formulas,
\begin{equation}\label{eq:13}cos(y)=\sqrt{\frac{1+cos(2y)}{2}} \equiv a_{n}=\sqrt{ \frac{1+a_{n-1}}{2}},
\end{equation}
and
\begin{equation}\label{eq:14}sin(y)=\sqrt{\frac{1-cos(2y)}{2}} \equiv c_{n}=\sqrt{ \frac{1-a_{n-1}}{2}},
\end{equation}
where $ n>1.$ \vspace{3mm}\newline
Given the first term  $a_1=\frac{s}{m},$
\begin{equation}\label{eq:15}
c_n=\left[\frac{1}{(2^{k_{n+1}}m)^{q_{n+1}}}
\underbrace{
\sqrt{ (2^{k_n}m)^{q_n}-\sqrt{(2^{k_{n-1}}m)^{q_{n-1}}+\cdots+\sqrt{m+s}}}}_{n \hspace{1mm}square \hspace{1mm} roots}\right],
\end{equation}
\begin{multline}\label{eq:16}
sin\left(\frac{\arccos\frac{s}{m}}{2^{n}}\right)=\\\left[\frac{1}{(2^{k_{n+1}}m)^{q_{n+1}}}
\underbrace{
\sqrt{ (2^{k_n}m)^{q_n}-\sqrt{(2^{k_{n-1}}m)^{q_{n-1}}+\cdots+\sqrt{m+s}}}}_{n \hspace{1mm}square \hspace{1mm} roots}\right]
\end{multline}
where $k_n=2^{n-1}-1$,  $q_n=2^{1-n},$    \vspace{3mm}\newline
Similarly to (\ref{eq:10}),
\begin{equation}\label{eq:17}
\lim_{n \to \infty} \left[\frac{2^{n}}{\arccos\frac{s}{m}}c_n\right]=1
\end{equation}
As required, (\ref{eq:17}) can be written as
\begin{equation}\label{eq:18}
\lim_{n \to \infty}\left[g\left(\frac{s}{m}, n\right)
\underbrace{
\sqrt{ (2^{k_n}m)^{q_n}-\sqrt{(2^{k_{n-1}})^{q_{n-1}}+\cdots+\sqrt{m+s}}}}_{n \hspace{1mm}square \hspace{1mm} roots} \right]=1,
\end{equation}
where,  $g\bigg(\frac{s}{m}, n\bigg)=\frac{2^{n}}{arcos\frac{s}{m}  \big((2^{k_{n+1}}m)^{q_{n+1}}\big)}$
\vspace{3mm}\newline
\textit{Examples} : The following are some examples of the general nested square root formula in (\ref{eq:16}), for $s=1,$ $m=5,7,9$ and $11$ respectively,
\begin{align}
1 =\lim_{n \to \infty} \left[ g\bigg(\frac{s}{m}, n\bigg)
\underbrace{
\sqrt{f(n,m)-\sqrt{f(n-1,m)+\cdots+\sqrt{40^{\frac{1}{4}}+\sqrt{10^{\frac{1}{2}}+\sqrt{6}}}}}}
_{n\hspace{1mm} square\hspace{1mm}roots}
\right]\label{eq:19}\\
1 =\lim_{n \to \infty} \left[ g\bigg(\frac{s}{m}, n\bigg)
\underbrace{
\sqrt{f(n,m)-\sqrt{f(n-1,m)+\cdots+\sqrt{56^{\frac{1}{4}}+\sqrt{14^{\frac{1}{2}}+\sqrt{8}}}}}}
_{n\hspace{1mm} square\hspace{1mm} roots}
\right]\label{eq:20},\\
1 =\lim_{n \to \infty} \left[ g\bigg(\frac{s}{m}, n\bigg)
\underbrace{
\sqrt{f(n,m)-\sqrt{f(n-1,m)+\cdots+\sqrt{72^{\frac{1}{4}}+\sqrt{18^{\frac{1}{2}}+\sqrt{10}}}}}}
_{n\hspace{1mm}square \hspace{1mm} roots}
\right]\label{eq:21},\\
1 =\lim_{n \to \infty} \left[ g\bigg(\frac{s}{m}, n\bigg)
\underbrace{
\sqrt{f(n,m)-\sqrt{f(n-1,m)+\cdots+\sqrt{88^{\frac{1}{4}}+\sqrt{22^{\frac{1}{2}}+\sqrt{12}}}}}}
_{n \hspace{1mm} square \hspace{1mm} roots}
\right]\label{eq:22},
\end{align}
where, $f(n,m)=(2^{k_{n}}m)^{q_{n}}$
\begin{remark}
It is clear that in examples  (\ref{eq:19}), (\ref{eq:20}), (\ref{eq:21}) and (\ref{eq:22})  the numbers within the nested square roots have fractional powers  as we sought to derive.
\end{remark}
\section{ Finitely Nested Square Roots}
In the previous section we derived nested square root formulas convergent to unity, by limiting the nth term of the recursion formula in (\ref{eq:14}) to infinity. In this section however, we discuss a  nested square root formulas convergent to unity for a finite nth term of the recursion formula in (\ref{eq:14}). To begin, let's establish the conditions  for the convergence. \newline
Let
\begin{equation}\label{eq:23}
 s^2=m^2-d^2,   
\end{equation}
where $d>0,$\newline \vspace{3mm}
Then, the first term of the recursion formula,
\begin{equation}\label{eq:24}
 a_1=\frac{s}{m}=\frac{\sqrt{m^2-d^2}}{m}=\sqrt{1-\frac{d^2}{m^2}}, \end{equation}
 Using Binomial expansion,
 \begin{equation}\label{eq:25}
   \left(1-\frac{d^2}{m^2}\right)^{\frac{1}{2}}=1-\frac{d^2}{2m^2}-\frac{d^4}{8m^4}-\frac{d^6}{16m^6}-\frac{5d^{8}}{128m^{8}}-\cdots 
 \end{equation}
 (\ref{eq:23}) converges to unity as $\frac{d}{m} \to  0.$ This condition can simply be achieved in two ways. We can keep $d$ constant while limiting m to infinity or keep $m$ constant while limiting $d$ to 0. In both cases  $\frac{\sqrt{m^2-d^2}}{m} \to 1,$ and $arccos(a_1)\to 0, $ 
\begin{equation}\label{eq:26}
\lim_{\substack{m \to \infty \\  d \to 0 }}\left[ g\bigg(\frac{s}{m}, n\bigg) c_n
\right]=1\end{equation}
\textit{Examples}: The following are some examples of finite nested square root formulas in (\ref{eq:26}), for $n=2,3$ and $4$ respectively,
\begin{equation}
  \lim_{\substack{m \to \infty \\  d \to 0 }}\left[g\left(\frac{s}{m}, 2\right) \sqrt{(2m)^{\frac{1}{2}}-\sqrt{m+s}}\right]=1,\end{equation}
  \begin{equation}
  \lim_{\substack{m \to \infty \\  d \to 0 }}\left[g\left(\frac{s}{m}, 3\right) \sqrt{(2^3m)^{\frac{1}{4}}-\sqrt{(2m)^{\frac{1}{2}}+\sqrt{m+s}}}\right]=1,
  \end{equation}
  \begin{equation}
  \lim_{\substack{m \to \infty \\  d \to 0 }}\left[g\left(\frac{s}{m}, 4\right) \sqrt{(2^7m)^{\frac{1}{8}}-\sqrt{(2^3m)^{\frac{1}{4}}+\sqrt{(2m)^{\frac{1}{2}}+\sqrt{m+s}}}}\right]=1,
\end{equation}
\section{Conclusion}
We have discussed the convergence of nested square root formulas in two ways-by infinite and finite recursion formulas.  We derived these formulas  separately.However by combining these methods, we derive a more efficient nested square root formula convergent to unity , in the form
\begin{equation}\label{eq:30}
\lim_{\substack{n \to \infty\\ m \to \infty \\ d \to 0}}\left[ g\left(\frac{s}{m}, n\right) c_n=1
\right],\end{equation}
with all the variables and functions as previously defined.\newline

\end{document}